\documentclass[12pt,a4paper]{amsart}
\usepackage{amsfonts,amsmath,amssymb}

\def\GL{\operatorname{\mathrm {GL}}\nolimits}
\def\SU{\operatorname{\mathrm {SU}}\nolimits}
\def\SO{\operatorname{\mathrm {SO}}\nolimits}
\def\SL{\operatorname{\mathrm {SL}}\nolimits}

\def\Lie{\operatorname{\mathrm {Lie}}\nolimits}
\def\Aut{\operatorname{\mathrm {Aut}}\nolimits}
\def\Gal{\operatorname{\mathrm {Gal}}\nolimits}
\def\Rad{\operatorname{\mathrm {Rad}}\nolimits}
\def\diag{\operatorname{\mathrm {diag}}\nolimits}
\def\sign{\operatorname{\mathrm {sign}}\nolimits}
\def\sgn{\operatorname{\mathrm {sign}}\nolimits}
\def\Cent{\operatorname{\mathrm {Cent}}\nolimits}
\def\vol{\operatorname{\mathrm {Vol}}\nolimits}
\def\Hei{\operatorname{\mathrm {Heis}}\nolimits}
\def\tr{\operatorname{\mathrm {trace}}\nolimits}
\def\ad{\operatorname{\mathrm {ad}}\nolimits}

\def\Id{\operatorname{\mathrm {Id}}\nolimits}
\def\ind{\operatorname{\mathrm {Ind}}\nolimits}
\def\bbR{\mathbb R}
\def\bbC{\mathbb C}
\def\bbS{\mathbb S}
\def\bbT{\mathbb T}
\newtheorem{thm}{Theorem}[section]
\newtheorem{prop}[thm]{Proposition}

\newtheorem{rem}[thm]{Remark}
\def\beginpf{\textsc{Proof.}}
\def\endpf{\hfill$\Box$}
\title[Poisson Summation and Endoscopy for $\SU(2,1)$]{Poisson Summation and Endoscopy\\ for $\SU(2,1)$}
\author[Do Ngoc Diep, Do Thi Phuong Quynh]{Do Ngoc Diep${}^{1}$, Do Thi Phuong Quynh${}^{2}$}
\address{${}^{1}$ Institute of Mathematics, VAST, 18 Hoang Quoc Viet road, Cau Giay district, 10307 Hanoi, Vietnam \newline
{\tt Email: dndiep@math.ac.vn}}
\address{${}^{2}$ {\sc Medicine University Colllege, Thai Nguyen Universtiy, Thai Nguyen City, Vietnam}\newline
{\tt Email:  phuongquynhtn@gmail.com}}

\begin{document}
\date{\bf Version of \today}

\maketitle
\begin{abstract} 
In this paper we analyze the endoscopy for $\SU(2,1)$. We make a precise realization of the discrete series representations (following the work \cite{liu} in Section 2), a computation of their traces (Section 3) and an exact formula for the  Poisson summation and endoscopy for this group (in Section 4).\\
\\
\textsl{Key words:} discrete series representation; Poisson summation; trace formula; endoscopy\\
\textbf{2010 Subject Classification: } Primary: 22E45, Secondary: 11F70; 11F72
\end{abstract}
\tableofcontents
\section{Introduction}

For a fixed reductive group  we do find, if possible all the irreducible unitary representations  and then decompose  a particular represention into a discrete direct sum or continuous sum (direct integral) of irreducible ones; they  are the basic questions of harmonic analysis on reductive Lie groups. In particlular, the discrete part of the regular representation of reductive groups is the discrete sum of discrete series representations. 

Often these problems are reduced to the trace formula, because as known, the unitary representations are uniquely defined by it generalized character and infinitesimal character. 
Following Harish-Chandra, the generalized character is defined by its restriction to the maximal compact subgroup, as the initial eigenvalue problem for the generalized Laplacian (the Casimir operator) with the infinitesimal character as the eigenvalues infinitesimal action of Casimir operators.
There is a very highly developed theory of Arthur-Selberg trace formula. The theory is complicated and one reduces it to the same problem for smaller   endoscopic subgroups. It is called the transfer and plays a very important role in the theory.
By definition, an endoscopic subgroup is the connected component of the centralizer of regular semisimple elements, associated to representations, namely by the orbit method.

In the previous papers \cite{diepquynh1},\cite{diepquynh2}, we treated the case of rank one group $\SL(2,\bbR)$ and the rank 2 group  $\SL(3,\bbR)$ of rank 2. In this paper we  do the same study for $\SU(2,1)$, also of rank 2.

For the discrete series representations of $\SU(2,1)$ the following endoscopic groups should be considered:
\begin{itemize}
\item \textit{the elliptic case}: diagonal subgroup of regular elements
$$H= \left\{\diag(a_1,a_2,a_3); \quad a_1a_2a_3 = 1, a_i\ne a_j,\forall i \ne j\right\}\cong \bbT^2;$$
\item \textit{the parabolic case}: blog-diagonal subgroup of regular elements 
$$H = \left\{\begin{pmatrix} \cos\theta & \sin\theta & 0\\ -\sin\theta & \cos\theta & 0 \\ 0 & 0 & 1\end{pmatrix}\right\}\cong \SO(2);$$
\item \textit{the trival case}: the group $H= \SU(2,1)$ its-self.
\end{itemize}
We show in this paper that the method that J.-P. Labesse \cite{labesse} used for $\SL(2,\bbR)$ is also applicable  for $\SU(2,1)$.
Therefore one deduces the transfer formula for the discrete series representations and discrete series limits of  $\SU(2,1)$ to the corresponding endoscopic groups. 

The content of the paper can be described as follows.
For the group $\SU(2,1)$ we make a precise realization of the discrete series representations (following \cite{liu} in Section 2) by using the Orbit Method and Geometric Quantization to the solvable radical, a computation in the context of $\SU(2,1)$ of their traces (Section 3) and an exact formula for the noncommutative Poisson summation and endoscopy of  for this group (in Section 4). 

\section{Irreducible Unitary Representations of $\SU(2,1)$}

\subsection{The structure of $\SU(2,1)$} The following notions and results are folklore and we recall them only to fix an appropriate system of notations.

Let us remind that the group $\SU(2,1)$ is 
$$\SU(2,1)= \left\{ X \in \GL(3,\mathbb C) | \overline{{}^tX} I_{2,1} X = I_{2,1}  \right\}$$ where $I_{2,1}$ is the matrix of the quadratic form
$$Q(u) = |u_1|^2 + |u_2|^2 - |u_3|^2, u \in \mathbb C^3,$$ i.e. 
$$I_{2,1} = \begin{pmatrix} 1 & 0 & 0\\ 0 & 1 & 0\\ 0 & 0 & -1 \end{pmatrix}.$$

Denote by $\mathfrak{su}(2,1)$ the Lie algebra $\Lie\SU(2,1)$, $\theta$ the Cartan involution of the group $G=\SU(2,1)$. The corresponding Cartan involution for its Lie algebra $\mathfrak{su}(2,1)$ is denote by the same symbol $\theta\in \Aut\mathfrak{su}(2,1)$, $\theta^2 = \Id$,  $$\theta(X) = \overline{{}^tX^{-1}}, X\in \SU(2,1).$$ 
The maximal compact subgroup $K$ of $G$
$$K = \left\{\left. \begin{pmatrix} x & 0\\ 0 & \det x^{-1}\end{pmatrix} \right| x \in U(2)    \right\}$$ 
 is the subgroup of $G$, the Lie algebra $\mathfrak k$ of which is consisting of all the matrices with eigenvalue $+1$ of the Cartan involution,
$$\mathfrak k = \left\{\left. \begin{pmatrix} A & 0\\ 0 & -\tr A\end{pmatrix} \right| A \in \mathfrak u(2)    \right\}.$$ 
The Borel subgroup of $\SU(2,1)$ is the minimal parabolic subgroup $$P_0 = B=\left\{\left. p= \begin{pmatrix} m & *\\ 0 & \det m^{-1} \end{pmatrix}\right| m \in U(2) \right\},$$  
the Lie algebra of which is consisting of all the matrices with eigenvalue $-1$ of the Cartan involution,,
$$\mathfrak b = \left\{\left.  X\in \mathfrak g \right| \theta(X) = -X   \right\}.$$

The group $\SU(2,1)$ admits the well-known Cartan decomposition in form of a  product
$G= BK$. The Borel subgroup  $B$ is endowed with a further decomposition into a semi-direct product $B=MA\ltimes U$ of a maximal split torus $A$, the Lie algebra of which is
$$\mathfrak a = \left\{\left. H = \begin{pmatrix} 0 & 0 & \lambda\\
0 & 0& 0\\ \lambda & 0 & 0\end{pmatrix} \right| \lambda \in \mathbb R  \right\}$$
 and the unipotent radical $U= \Rad_uB \cong \Hei_3$, generated by matrices
$$X = \begin{pmatrix} 
0 & -1 & 0\\ 
1 & 0 & -1\\
0 & -1 & 0\end{pmatrix},
Y = \begin{pmatrix} 
0 & -1 & 0\\ 
-1 & 0 & 1\\
0 & -1 & 0\end{pmatrix},
Z = \begin{pmatrix} 
2 & 0 & -2\\ 
0 & 0 & 0\\
2 & 0 & -2\end{pmatrix},$$
satisfying the Heisenberg commutation relation $[X,Y]=Z, [X,Z] = [Y,Z] = 0$.
$$\mathfrak b = \mathfrak m \oplus \mathfrak a \oplus \mathfrak u,$$ where 
$\mathfrak u =\Lie \Hei_3 = \langle X,Y,Z\rangle$, $\mathfrak a = \langle H \rangle$, $\mathfrak m = \langle T=\frac{1}{3}\diag(1,-2,1)\rangle$. We have the commutative relations of harmonic oscillator $$[T,X] = Y, [T,Y]=-X, [T,Z] = 0$$

The center $C(\mathfrak k)$ is of dimension 1,
$$C(\mathfrak k) = \left\{\left. \begin{pmatrix}
i\lambda & 0 & 0\\
0 & i\lambda & 0\\
0 & 0 & -2i\lambda 
\end{pmatrix}  \right| \lambda \in \mathbb R, i = \sqrt{-1}  \right\}$$
There is a compact Cartan subalgebra $\mathfrak h$ consisting of all diagonal matrices $$\mathfrak h = \left\{\left.  \diag(ih_1,ih_2,ih_3 ) \right|  h_1,h_2,h_3 \in \mathbb R, h_1 + h_2 + h_3 = 0\right\}\subset \mathfrak k.$$
The associate root system is 
$$\Delta(\mathfrak g_\mathbb C , \mathfrak k_\mathbb C) = \left\{\left. \alpha_{kl} = \alpha_k - \alpha_l \right| \alpha_k(h_l) = \delta_{kl}\right\}.$$
It means that $\alpha_{kl} = (0,\dots,0,\underbrace{1}_k,0,\dots, 0,\underbrace{-1}_l,0,\dots,0 )\in \mathfrak h^* $, $1 \leq k \ne l \leq 3$.
The subroot system of compact roots is
$\Delta_c = \{\pm \alpha_{12}\} $ The noncompact root system is
$\Delta_n = \{\pm\beta,\pm 2\beta\}$  where $\beta$ is the noncompact root such that 
$$\beta( \begin{pmatrix} 0 & 0 & \lambda\\
0 & 0& 0\\ \lambda & 0 & 0\end{pmatrix}) = \lambda$$ 

The coroot system is 
$$\Delta(\mathfrak g_\mathbb C , \mathfrak k_\mathbb C)^* = \left\{ H_{kl} = E_{kk} - E_{ll} \left| {E_{kl} = \mbox{ elementary matrix with the }\atop \mbox{only nonzero entry 1 on position } (i,j)}\right. \right\}$$
 
\begin{prop}[Orbit picture]
The space $\mathfrak b^*$ is divided into a dijoint union of the following coadjoint $B$-orbits:
\begin{enumerate}
\item[a.] The two half-spaces $\Omega_\pm$ consisting of functionals $F = tT^* +  xX^* + yY^* +zZ^*$, 
$$\Omega_+= \{(t,x,y,z) \in \mathbb R^4 |  z >  0\}$$ 
$$\Omega_-= \{(t,x,y,z) \in \mathbb R^4 |  z <  0\}$$ 
\item[b.] A family of cylinders with hyperbolic base
$$\Omega_\alpha = \{(t,x,y,0) \in \mathbb R^3 \times \{0\} | xy = \alpha  \}, \alpha > 0$$ 
\item[c.] Four half-planes corresponding to the case $xy = 0$ but $x\ne y$
$$\Omega_{x>0} = \{t,x,0,0) \in \mathbb R^4 | x > 0\}$$
$$\Omega_{x<0} = \{t,x,0,0) \in \mathbb R^4 | x < 0\}$$
$$\Omega_{y>0} = \{t,0,y,0) \in \mathbb R^4 | y> 0\}$$
$$\Omega_{y<0} = \{t,0,y,0) \in \mathbb R^4 | y< 0\}$$
\item[d.]
The origin $$\Omega = \{(0,0,0,0) \}$$
\end{enumerate}
\end{prop}
\textsc{Proof.}
This proposition is proven by a direct computation.
\hfill$\Box$

Consider the linear functionals $\pm Z^*\in \Omega_\pm \subset \mathfrak g^*$ and the corresponding coadjoint orbits $\Omega_\pm = G.(\pm Z^*).$

\begin{prop}
Subalgebras $\mathfrak l = \mathbb C(X \pm iY) \oplus \mathbb CZ \subset \mathfrak u_\mathbb C$ are the positive polarizations at $\pm Z^*\in \Omega_\pm$.
\end{prop}
\textsc{Proof.}
The proposition is proven by a direct computation.
\hfill$\Box$
\begin{rem}
The coadjoint action of $K$ in $\mathfrak g^*$ keeps the set $\mathfrak b^*$ invariant.
\end{rem}

\subsection{Holomorphic Induction}
Following the orbit method and the holomorphic induction, we do choose the integral functionals $\lambda$, take the corresponding orbits and then choose polarization and use the holomorphic induction.

As described above, the positive root system $\Delta^+ = \Delta^+_c \cup \Delta^+_n = \{ \alpha_{kl}, 1\leq k\ne l \leq 3, \beta, 2\beta\}=\{\alpha_{12}, \alpha_{32},\alpha_{31}\}, \rho = \alpha_{32}$.
The root spaces are $\mathfrak g_\mathbb C^{\alpha_{kl}} = \mathbb C   E_{kl}$. $\mathfrak g_\beta = \mathbb R X \oplus \mathbb R Y$ and $\mathfrak g_{2\beta} = \mathbb R Z$.  
Define $$\mathfrak p_+ = \bigoplus_{\alpha \in \Delta^+_n} \mathfrak g^\alpha = \mathfrak g^{\alpha_{32}} \oplus \mathfrak g^{\alpha_{31}}=\mathbb C E_{31}\oplus \mathbb C E_{32}$$ 
and 
$$\mathfrak p_- = \bigoplus_{\alpha \in \Delta^-_n} \mathfrak g^{\alpha} = \mathfrak g^{\alpha_{23}} \oplus \mathfrak g^{\alpha_{13}}=\mathbb C E_{13}\oplus \mathbb C E_{23}$$ 

Denote $\mathcal F\subset (i\mathfrak h)^*$ the set of all linear functional $\lambda$ on $\mathfrak h_\mathbb C$ such that $(\lambda + \rho)(H_\alpha)$ is integral for any root $\alpha\in\Delta$, where $H_\alpha$ is the coroot corresponding to root $\alpha$ and $\rho$ is the half-sum of the positive roots. 
Denote also
$$\mathfrak F' = \{\lambda \in \mathfrak F | \lambda(H_\alpha) \ne 0, \forall \alpha\in \Delta \},$$
$$\mathfrak F'_0 = \{\lambda\in \mathfrak F' | \lambda(H_\alpha) > 0, \forall \alpha\in \Delta^+ _c\}$$ 
$$=\left\{\lambda\in \mathfrak F' \left| \begin{array}{l} \lambda(H_{12}) \in \mathbb N^+ \mbox{ and } \lambda(H_{31}) \in \mathbb N^+  \mbox{ (holomorphic case) or }\\
\lambda(H_{12})\in \mathbb N^+ \mbox{ and } \lambda(H_{23}) \in \mathbb N^+  \mbox{(anti-holomorphic case)  or }\\
\lambda(H_{12}) \in \mathbb N^+ \mbox{ and } \lambda(H_{13}) \in \mathbb N^+, \lambda(H_{12}) > \lambda(H_{13}) \mbox{ neither-nor case }
\end{array}
 \right .\right\}$$
 
Choose complex subalgebra $$\mathfrak e = \mathfrak p_+ \oplus \mathfrak k_\mathbb C,$$ we have
$$\mathfrak e + \overline{\mathfrak e} = \mathfrak g_\mathbb C, \mathfrak e \cap \overline{\mathfrak e} = \mathfrak k_\mathbb C $$
and therefore we have a positive polarization.

The compact root Weyl group $W_K = \langle s_{\alpha_{12}} \rangle$ is generated by a single reflection $s_{\alpha_{12}}$ then
for any $\lambda \in i\mathfrak h$, $-s_{\alpha_{12}}\lambda - \alpha_{32} = -s_{\alpha_{12}}(\lambda + \alpha_{31})$, therefore if $V_\lambda$ is a $K$-module of lowest weight $\lambda + \rho$ then its contragradient $K$-module $V^*_\lambda$ is of heighest weight $\lambda + \alpha_{31}$.

Because $G= BK = B_1K$, the relative cohomology of  $(\mathfrak g, K)$-module with coefficients in the representation $V_\lambda$ can be reduced to the one of $B$ or $B_1=AU \subset B$ with Lie algebra $\mathfrak b_1 = \langle S= E_{13}+E_{31},X,Y,Z\rangle$. 
\begin{prop} For the $(\mathfrak g, K)$-module $\pi_\lambda$ with coefficients in the representation $V_\lambda$
$$H^{q_\lambda}(G,K;\mathfrak e,V_\lambda) = H^{q_\lambda}(B, M;\mathfrak e \cap \mathfrak b,V_\lambda) =H^{q_\lambda}(B_1;\mathfrak e \cap \mathfrak b_1,V_\lambda)$$
\end{prop}

\subsection{Hochschild-Serre spectral sequence}
Remark that because in general $\mathfrak p$ is not a subalgebra, we can modify it by taking subalgebra $\mathfrak h_+ = \mathbb C(Y+iX) \oplus \mathbb C(S-iZ/2)$:
$$\mathfrak e = \mathfrak p_+ \oplus \mathfrak k_\mathbb C = \mathfrak h_+ \oplus \mathfrak k_\mathbb C, \quad\mathfrak h_+=\mathfrak p_+ \cap \mathfrak h, \quad \mathfrak b_1 := \mathfrak b \cap \mathfrak p_+.$$ Therefore, one has
$$\mathfrak e \cap \mathfrak b_1 = \mathfrak h_+, \quad \mathfrak e \cap \mathfrak b = \mathfrak h_+ \oplus \mathfrak m_\mathbb C.$$
We have then a filtration on polarizable parabolic subalgebra. 
Therefore, we may construct a Hochschild-Serre spectral sequence for this filtration.

Consider a highest weight $\lambda + \alpha_{31}$ representation  $V^*_\lambda$ of $\mathfrak k_\mathbb C$, which is trivially on $\mathfrak p_+$ extended to a representation $\xi$ of $\mathfrak e = \mathfrak p_+ \oplus \mathfrak k_\mathbb C$. The action of $\mathfrak h_+$ in $V^{\lambda + \alpha_{31}}$ is $\xi + \frac{1}{2}\tr\ad_{\mathfrak b_1}$.
Denote by $\mathcal H_\pm$ the space of representations $T_\pm$ of $B$ $\Omega_\pm$ above and by $\mathcal H_\pm^\infty$ the subspaces of smooth vectors.
Because $\dim_\mathbb C(\mathfrak p_\mathbb C) = 2$, we have $\wedge^q(\mathfrak h_+)= 0$, for all $q \geq 3$.
It is natural to define the Hochschild-Serre cobound operators
$$(\delta_\pm)_{\lambda,q} : \wedge^q(\mathfrak h_+)^* \otimes V^{\lambda + \alpha_{31}} \otimes \mathcal H_\pm^\infty \to  \wedge^{q+1}(\mathfrak h_+)^* \otimes V^{\lambda + \alpha_{31}} \otimes \mathcal H_\pm^\infty$$ and by duality their formal adjoint operators $(\delta_\pm)_{\lambda,q}^*$.
The Hochschild-Serre spectral sequence is convergent
$$\bigoplus_{r+s=q} H^r(\mathfrak e_1; H^s(M;V^{\lambda + \alpha_{31}} \otimes \mathcal H_\pm^\infty)) \Longrightarrow H^q(B;\mathfrak b_1,V_\lambda)$$

\subsection{Holomorphic or anti-holomorphic case} 
Denote by $\mathfrak H_{\lambda,\pm}^q$ the intersection $\ker(\delta_\pm)_{\lambda,q}) \cap \ker(\delta_\pm)_{\lambda,q})^*$ 
The main result  from \cite{liu}, what we do not need to use  in the rest of the paper,  is
\begin{thm}[\cite{liu}]
a. $\pi_\lambda|_{B_1} = \pi^{q_\lambda}(\mathfrak e)|_{B_1} = \dim (\mathfrak H_{\lambda,-}^{q_\lambda})T_+ \oplus \dim (\mathfrak H_{\lambda,+}^{q_\lambda})T_-$

b. $\pi_\lambda|_{B} = \pi^{q_\lambda}(\mathfrak e)|_{B} = \sum_{m\in\mathbb Z}[\dim (\mathfrak H_{\lambda,-}^{q_\lambda})_{\sigma_m}^{q_\lambda}T_{m,+} \oplus \dim (\mathfrak H_{\lambda,+}^{q_\lambda})_{\sigma_m}^{q_\lambda}T_{m,-}]$
\end{thm}

\subsection{Neither holomorphic nor anti-holomorphic case}
\begin{thm}[\cite{liu}]
If $\pi_\lambda$ is in the neither holomorphic nor anti-holomorphic cases, we have 
$$\pi_\lambda|_B = \sum_{m=0}^\infty  T_{[3m + \frac{3f_\lambda(H) - f_\lambda(Z)}{2}],+} \oplus  
\sum_{m=0}^\infty  T_{[-(3m + \frac{-(3f_\lambda(H) + f_\lambda(Z)}{2})],-}$$ 
\end{thm}

\section{Trace Formula}
The main idea to compute the traces as characters of the discrete series representation is the following result of I. M. Gelfand and M.A. Naimark
\begin{thm}[I. M. Gelfand, M. A. Naimark]
For any irreducible unitary representation $\pi$ of reductive group $G$, the Gelfand image operator $\pi(f)$ is of trace class for any function $f$ with compact support on $G$.
\end{thm}
Because of the properties of trace class operators, the character $\chi_\pi(f)$ is a distribution of $G$ and therefore, following I.M. Gelfand is called the character of representation.

The discrete part of the regular representation is decomposed into a discrete orthogonal sum of irreducible representations and there the trace is reduce to a sum of all characters of the discrete series characters $\Theta_{n,\pm}(f)$.
$$\tr {}^oR(f)= \sum_{{\pi\in \mbox{ disrete series }\atop \varepsilon = \pm 1}} \sum \Theta_{n,\varepsilon}(f)$$

\subsection{Trace formula}
Let us remind that $\Gamma \subset \SU(2,1)$ is a finitely generated Langlands type discrete subgroup with finite number of cusps.
Let $f\in C^\infty_c(\SU(2,1))$ be a smooth function of compact support. If $\varphi$ is a function from the representation space, the action of the induced representation $\ind_P^G\chi$ is the restriction of the right regular representation $R$ on the  space of induced representation.
$$ R(f)\varphi = \int_G (f(y)R(y)\varphi(x))dy = \int_G f(y)\varphi(xy)dy$$
$$= \int_G f(x^{-1}y)\varphi(y)dy (\mbox{right invariance of Haar measure } dy) $$
$$= \int_{\Gamma\backslash G}\left(\sum_{\gamma\in \Gamma}f(x^{-1}\gamma y) 
\right)\varphi(y)dy$$

 Therefore, this action can be represented by an operator with kernel $K(x,y)$ of form
$$[R(f)\varphi](x) = \int_{\Gamma\backslash G} K_f(x,y)\varphi(y)dy,$$ where
$$K_f(x,y) = \sum_{\gamma\in \Gamma}f(x^{-1}\gamma y).$$ Because the function $f$ is of compact support, this sum is convergent, and indeed is a finite sum, for any fixed $x$ and $y$ and is of class $L^2(\Gamma\backslash G \times\Gamma\backslash G )$. The operator is of trace class and it is well-known that
$$\tr R(f) = \int_{\Gamma\backslash G}K_f(x,x)dx.$$
As supposed, the discrete subgroup $\Gamma$ is finitely generated. Denote by $\{\Gamma\}$ the set of representatives of conjugacy classes. For any $\gamma\in\Gamma$ denote the centralizer of $\gamma\in \Omega \subset G$ by $\Omega_\gamma$, in particular, $G_\gamma \subset G$. Following the Fubini theorem for the double integral, we can change the order of integration to have
$$\tr R(f) = \int_{\Gamma\backslash G} K_f(x,x)dx =\int_{\Gamma\backslash G} \sum_{\gamma\in \Gamma} f(x^{-1}\gamma x)dx$$
$$=\int_{\Gamma\backslash G} \sum_{\gamma\in \{\Gamma\}}\sum_{\delta\in \Gamma_\gamma\backslash \Gamma} f(x^{-1}\delta^{-1}\gamma\delta x)dx$$
$$=\sum_{\gamma\in \{\Gamma\}}\int_{\Gamma_\gamma\backslash G}  f(x^{-1}\gamma x)dx=\sum_{\gamma\in \{\Gamma\}}\int_{G_\gamma\backslash G}\int_{\Gamma_\gamma\backslash G_\gamma}  f(x^{-1}u^{-1}\gamma u x)dudx$$
$$= \sum_{\gamma\in \{\Gamma\}} \vol(\Gamma_\gamma \backslash G_\gamma)\int_{G_\gamma\backslash G}f(x^{-1}\gamma x)dx.$$
Therefore, in order to compute the trace formula, one needs to do:
\begin{itemize}
\item classfiy the conjugacy classes of all $\gamma$ in $\Gamma$: they are of type elliptic (different eigenvalues of the same sign),
 hyperbolic (nondegenerate, with eigenvalues of different sign),
 parabolic (denegerate)  
\item Compute the volume of form $\vol(\Gamma_\gamma \backslash G_\gamma)$; it is the volume of the quotient of the stabilazer of the adjoint orbits. 

\item and compute the orbital integrals of form
$$\mathcal{O}_\gamma(f) = \int_{G_\gamma\backslash G} f(x^{-1}\gamma x)d\dot x$$ 
The idea is to reduce these integrals to smaller endoscopic subgroups in order to the correponding integrals are ordinary or almost  ordinary. 
\end{itemize}

\subsection{Stable trace formula}
The Galois group $\Gal(\mathbb C/\mathbb R) \cong \mathbb Z_2$ of the complex field $\mathbb C$ is acting on the discrete series representation by character $\kappa(\sigma) = \pm 1$. The change of sign in $\mathfrak A^*$ changes the orientation of the space $\mathfrak A^* \cong \mathbb R^2$.  Therefore the sum of characters can be rewrite as some sum over stable classes of characters. Following Harish-chandra parametrization of the discrete series representations $\pi_{w\lambda + \rho}$, where $w\in W$, $\lambda\in C=\mathfrak A^*_+$ the positive chamber in $\mathfrak A^*$,
$$\tr R(f) = \sum_{w,\lambda} \sum_{\varepsilon = \pm 1} (\Theta_{w\lambda + \rho,+}(f) - \Theta_{w\lambda + \rho,-}(f)).$$

\section{Endoscopy and Poisson Summation}
In order to compute the orbital integral
$$\mathcal O_\gamma(f) = \int_{G_\gamma\backslash G} f(x^{-1}\gamma x)dx$$ we do: 
\begin{itemize}
\item classify all possible cases of $\gamma$: all distinguished eigenvalues, satisfying the condition $a_1a_2a_3 = 1$ (the simplest case), and two of eigenvalues are concided $a_1=a_2$ and $a_1^2a_3= 1$ (the second case).
\item  transform integration to reduce to the corresponding possible endoscopic subgroups.
\end{itemize} 

We observe that the method, J.-P. Labesse \cite{labesse} proposed for $\SL(2,\bbR)$ is still applicable for $\SU(2,1)$ with some modifications.

\subsection{Orbital integrals}
\textsl{The simplest case} is the elliptic case when $\gamma = \diag(a_a,a_2 ,a_3), \quad a_1a_2a_3 = 1$ with pairwise distinguished $a_i$. In this case, because of Iwasawa decomposition $x=mauk$, and the $K$-bivariance,  the orbital integral is
$$\mathcal O_\gamma(f) = \int_{G_\gamma\backslash G} f(x^{-1}\gamma x)dx = \int_U f(u^{-1}\gamma u)du$$ $$=\int_{\mathbb R^3} f(\begin{pmatrix} 1 & x & z\\ 0 & 1& y \\ 0 & 0 &  1\end{pmatrix}^{-1}\begin{pmatrix} a_1 & 0 & 0\\ 0 & a_2 & 0 \\ 0 & 0 & a_3 \end{pmatrix} \begin{pmatrix} 1 & x & z\\ 0 & 1& y \\ 0 & 0 &  1\end{pmatrix})dxdydz$$ 
$$=\int_{\mathbb R^3} f(\begin{pmatrix} 1 & -x & yx - z\\ 0 & 1 & -y\\ 0 & 0 & 1\end{pmatrix}\begin{pmatrix} a_1 & 0 & 0\\ 0 & a_2 & 0\\ 0 & 0 & a_3 \end{pmatrix} \begin{pmatrix} 1 & x & z\\ 0 & 1& y\\ 0 & 0 & 1\end{pmatrix})dxdydz$$ $$= \int_{\bbR^3} f(\begin{pmatrix} a_1  &(a_1-a_2)x & (a_1-a_3)z+a_3xy\\ 0 &a_2 & (a_2-a_3)y\\ 0 & 0 & a_3\end{pmatrix})dxdydz $$ $$= |a_1-a_2|^{-1}|a_2 -a_3|^{-1}|a_2 - a_3|^{-1} \mathcal O_\gamma(\tilde{f}).$$
The integral is abosolutely and uniformly convergent and therefore is a smooth function of $a\in \bbR^*_+$. Therefore the function $$f^H(\gamma) = \Delta(\gamma)\mathcal O_\gamma(f), \quad \Delta(\gamma) =\prod_{1\leq i<j\leq 3} |a_i- a_j|$$ is a smooth function on the endoscopic group $H=A= (\bbR^*)^2$. 

\textsl{The second case} is the case where $\gamma = \begin{pmatrix} k_\theta & 0\\ 0 & 1\end{pmatrix} = \begin{pmatrix} \cos\theta & \sin\theta & 0\\ -\sin\theta & \cos\theta & 0\\ 0 & 0 &1 \end{pmatrix}$. We have again, $x= mauk$ and
$$\mathcal O_{k(\theta)}(f) = \int_{G_{k(\theta)}\backslash G} f(k^{-1}u^{-1}a^{-1}m^{-1}k(\theta)mauk)dmdudadk$$
$$= \int_{G_{k(\theta)}\backslash G} f(u^{-1}a^{-1}m^{-1}k(\theta)mau)dmduda$$
$$ = \int_{G_{k(\theta)}\backslash G} f(\begin{pmatrix} 1 & -x &yx-z\\ 0 & 1 & -y\\ 0 & 0 & 1 \end{pmatrix}\begin{pmatrix} a_1^{-1} & 0 & 0\\ 0 & a_1^{-1} & 0\\ 0 & 0 & a_3^{-1} \end{pmatrix} \begin{pmatrix} \cos\theta & \sin\theta & 0\\ -\sin\theta & \cos\theta & 0\\ 0 & 0 & 1 \end{pmatrix}$$ $$\times \begin{pmatrix} a_1 &0 & 0\\ 0 & a_1 & 0 \\ 0 & 0 & a_3 \end{pmatrix}\begin{pmatrix} 1 & x & z \\ 0 & 1 & y \\ 0 & 0 & 1  \end{pmatrix})duda $$
$$= \int_1^\infty \tilde{f}(\begin{pmatrix}\cos\theta & t_1\sin\theta & 0\\ -t_1^{-1}\sin\theta & \cos\theta & 0\\ 0 & 0 & 1\end{pmatrix})\prod_{i=1}^3|t_i-t_i^{-1}|\frac{dt_i}{t_1}$$ 
$$= \int_0^{+\infty}\sgn(t-1) \tilde{f}(\begin{pmatrix}\cos\theta & t\sin\theta & 0\\ -t^{-1}\sin\theta & \cos\theta & 0\\ 0 & 0 & 1 \end{pmatrix})dt= F(\sin\theta).$$ 
When $f$ is an element of the Hecke algebra, i.e. $f$ is of class $C^\infty_0(G)$ and is $K$-bivariant, the integral is converging absolutely and uniformly. Therefore the result is a function $F(\sin\theta)$. 
The function $f$ has compact support, then the integral is well convergent at $+\infty$. At the another point $0$, we develope the function $f$ into the Taylor-Lagrange of the first order with respect to $\lambda = \sin\theta \to 0$, $\tilde{f}(t) = \tilde{f}(0) + t\tilde{f'}(0) + tB_{\tilde{f}}(t)$, and therefore after intergrating we have
$$F(\lambda) = A(\lambda) + \lambda B(\lambda),$$ where $A(\lambda) = f(0)$ and $B(\lambda)$ is the error-correction term $F'(\tau)$ at some intermediate value $\tau, 0 \leq \tau \leq t$.
Remark that 
$$ \begin{pmatrix} \sqrt{1-\lambda^2} & t\lambda & 0\\ -t^{-1}\lambda & \sqrt{1-\lambda^2}& 0 \\ 0 & 0 & 1 \end{pmatrix}$$ 
$$=\begin{pmatrix} t^{1/2}&0& 0\\ 0 & t^{-1/2} & 0\\
 0 & 0 & 1 \end{pmatrix} \begin{pmatrix} \sqrt{1-\lambda^2} & \lambda & 0\\ -\lambda & \sqrt{1-\lambda^2} & 0\\ 0 & 0 & 1 \end{pmatrix}\begin{pmatrix} t^{-1/2} & 0 & 0\\ 0 &t^{1/2} & 0\\
0 & 0 & 1 \end{pmatrix} $$ we have
$$B= \frac{dF(\tau)}{d\lambda} =\frac{d}{d\lambda} \left. \int_0^{+\infty} \sgn(t-1)f(\begin{pmatrix} \sqrt{1-\lambda^2} & t\lambda & 0\\ -t^{-1}\lambda & \sqrt{1-\lambda^2} & 0\\  0 & 0 & 1 \end{pmatrix})dt \right|_{t=\tau}$$
$$=  \int_0^{+\infty} \sgn(t-1)g(\begin{pmatrix} \sqrt{1-\lambda^2} & t\lambda & 0\\ -t^{-1}\lambda & \sqrt{1-\lambda^2}& 0\\ 0 & 0 & 1 \end{pmatrix})\frac{dt}{t}, $$
where $g\in C^\infty_c(N)$
and $g(\lambda) \cong O(-t^{-1}\lambda)^{-1}.$ $B$ is of logarithmic growth and $$B(\lambda) \cong \ln(|\lambda|^{-1})g(1)$$ up to constant term,  and therefore is contimuous. 
$$A = F(0) = |\lambda|^{-1} \int_0^\infty f(\begin{pmatrix}1 & \sgn(\lambda)u & 0\\ 0 & 1 & 0\\ 0 & 0 & 1\end{pmatrix}) du - 2f(I_3) + o(\lambda)$$
Hence the functions 
$$G(\lambda) = |\lambda|(F(\lambda) + F(\lambda)),$$
$$H(\lambda) = \lambda(F(\lambda) - F(-\lambda))$$ have the Fourier decomposition
$$G(\lambda) = \sum_{n=0}^N (a_n|\lambda|^{-1} + b_n)\lambda^{2n} + o(\lambda^{2N})$$
$$H(\lambda) = \sum_{n=0}^N h_n\lambda^{2n} + o(\lambda^{2N})$$
Summarizing the discussion, we have that in the case of $\gamma = k(\theta)$, there exists also a continuous function $f^H$ such that $$f^H(\gamma) = \Delta(\gamma) (\mathcal O_\gamma(f) - \mathcal O_{w\gamma}(f))=\Delta(k(\theta))\mathcal {SO}_\gamma(f), $$ where $\Delta(k(\theta)) = -2i\sin\theta$.

\subsection{Stable orbital integral}
Let us remind that the \textit{orbital integral} is defined as
$$\mathcal{O}(f) = \int_{G_\gamma\backslash G} f(x^{-1}\gamma x)d\dot x$$ 

The complex Weyl group is isomorphic to $\mathfrak S_3$ while the real Weyl group
is isomorphic to $\mathfrak S_2$ . The set of conjugacy classes inside a strongly regular
stable elliptic conjugacy class is in bijection with the pointed set
$\mathfrak S_3 /\mathfrak S_2$
 that can be viewed as a sub-pointed-set of the group
$\mathfrak E(\mathbb R, T, G) = (Z_2 )^2$
We shall denote by $\mathfrak K(\mathbb R, T, G)$ its Pontryagin dual.

Consider $\kappa \ne 1$ in $\mathfrak K(\mathbb R, T, G)$ such that $\kappa(H_{13}) = 1$. 
 Such a $\kappa$ is unique:
in fact one has necessarily
$\kappa(H_{12}) = \kappa(H_{13}) = -1$.

The endoscopic group $H$ one associates to $\kappa$ is isomorphic to
$S(U (1, 1) \times U (1))$
the positive root of $\mathfrak h$ in $H$ (for a compatible order) being $\alpha_{23} = \rho$

The endoscopic group $H$ can be embedded in $G$ as 
$$g(u,v,w) = \begin{pmatrix}
ua & iub & 0\\  -iuc & ud & 0\\  0 & 0 & v \end{pmatrix}, w = \begin{pmatrix} a & b\\ c & d \end{pmatrix}, ad-bc = 1 \mbox{ and } |u| = |v| = 1, uv = 1.$$
It will be useful to consider also the two-fold cover
$H_1 = S(U (1) \times U (1)) \times \SL(2)$.

Let $f_\mu$ be a pseudo-coefficient for the discrete series representation $\pi_\mu$ then the \textit{$\kappa$-orbital integral} of a  regular element $\gamma$
in $T (R)$ is given by
$$\mathcal O^\kappa_\gamma(f_\mu) =\int_{G_\gamma\backslash G} \kappa(x) f_\mu(x^{-1}\gamma x) d\dot x = \sum_{\sign(w) =1} \kappa(w) \Theta^G_\mu(\gamma^{-1}_w)  = \sum_{\sign(w) =1} \kappa(w)\Theta_{w\mu}(\gamma^{-1}),$$ because there is a natural bijection between the left coset classes and the right coset classes.

\begin{thm}\cite{labesse}
There is a natural function $\varepsilon : \Pi \to \pm 1$ such that in the Grothendieck group of discrete series representation ring, $$\sigma_G = \sum_{\pi\in \Pi} \varepsilon(\pi)\pi,$$ the map $\sigma \mapsto \sigma_G$ is dual to the map of geometric transfer, that for any $f$  on $G$, there is a unique $f^H$ on $H$
$$\tr \sigma_G(f) = \tr\sigma(f^H).$$
\end{thm}
\textsc{Proof.}
There is a natural bijection $\Pi_\mu \cong \mathfrak D(\mathbb R,H,G)$, we get a pairing
$$\langle .,.\rangle : \Pi_\mu \times \mathfrak K(\mathbb R,H,G)\to \mathbb C.$$ Therefore we have
$$\tr \Sigma_\nu(f^H) = \sum_{\pi\in \Pi_\Sigma} \langle s,\pi\rangle \tr \pi(f).$$
\hfill$\Box$

\subsection{Endoscopic transfer}
The transfer factor $\Delta(\gamma, \gamma_H )$ is given by
$$\Delta(\gamma, \gamma_H )=(−1)^{q(G)+q(H)} \chi_{G,H} (\gamma)\Delta_B (\gamma^{−1} ) . \Delta_{B_H}(\gamma_H^{-1})^{-1}$$
for some character $\chi_{G,H}$ defined as follows. Let $\xi$ be a character of the twofold covering $\mathfrak h_1$ of $\mathfrak h$, then
$\chi_{G,H} (\gamma^{−1} ) = e^{\gamma^{\rho - \rho_H +\xi}}$  defines a character of $H$, corresponding to $\mathfrak h$, because it is trivial on any fiber of the cover.

With such a choice we get when $\sign (w) = 1$
and $w\ne 1$, we have $\kappa(w) = -1$ and 
$$\Delta(\gamma^{-1},\gamma_H^{-1})\Theta_{w\mu}^G(\gamma) = -\frac{\gamma_H^{w\mu + \xi} - \gamma_H^{w_0w\mu + \xi}}{\gamma^{\rho_H}\Delta_{B_H}(\gamma_H}$$ therefore
$$\Delta(\gamma,\gamma_H)\Theta^G_{w\mu} (\gamma^{-1}) = \kappa(w)^{-1}\mathcal{SO}_\nu^H(\gamma_H^{-1}),$$ where $\nu = w\mu +\xi$ is running over the corresponding $L$-package of discrete series representations for the endoscopic group $H$.
Therefore we have the following formula
$$\Delta(\gamma,\gamma_H){\mathcal O}_\gamma^\kappa(f_\mu) = \sum_{\nu=w\mu + \xi\atop \sign(w) = 1} {\mathcal {SO}}_\nu^H(\gamma_H^{-1})$$ or
$$\Delta(\gamma,\gamma_H){\mathcal O}_\gamma^\kappa(f_\mu) = \sum_{\nu=w\mu + \xi\atop \sign(w) = 1} {\mathcal {SO}}_{\gamma_H}(g_\nu), $$ where $g_\nu$ is pseudo-coefficient for any one of the discrete series representation of the endoscopic subgroup $H$ in the $L$-package of $\mu$.

For any $$f^H = \sum_{\nu=w\mu+\rho\atop \sign(w) = 1} a(w,\nu)g_\nu,\quad a(w_1,w_2\mu = \kappa(w_2) \kappa(w_2w_1)^{-1}$$ we have the formula
$$\tr \Sigma_\nu (f^H) = \sum_{w} a(w,\nu) \tr \pi_{w\mu}(f).$$

\subsection{Poisson Summation and Endoscopy}
\begin{thm}\cite{labesse}
There is a natural function $\varepsilon : \Pi \to \pm 1$ such that in the Grothendieck group of discrete series representation ring, $$\sigma_G = \sum_{\pi\in \Pi} \varepsilon(\pi)\pi,$$ the map $\sigma \mapsto \sigma_G$ is dual to the map of geometric transfer, that for any $f$  on $G$, there is a unique $f^H$ on $H$
$$\tr \sigma_G(f) = \tr\sigma(f^H).$$
\end{thm}
\textsc{Proof.}
There is a natural bijection $\Pi_\mu \cong \mathfrak D(\mathbb R,H,G)$, we get a pairing
$$\langle .,.\rangle : \Pi_\mu \times \mathfrak K(\mathbb R,H,G)\to \mathbb C.$$ Therefore we have
$$\tr \Sigma_\nu(f^H) = \sum_{\pi\in \Pi_\Sigma} \langle s,\pi\rangle \tr \pi(f).$$
\hfill$\Box$

Suppose given a complete set of endoscopic groups $H = \mathbb S^1 \times \mathbb S^1 \times \{\pm 1\}$ or $\SL(2,\mathbb R) \times \{\pm 1\}$.  For each group, there is a natural inclusion
$$\eta: {}^LH \hookrightarrow {}^LG$$

Let $\varphi: DW_\bbR  \to {}^LG$ be the Langlands parameter, i.e. a homomorphism from the Weil-Deligne group $DW_\bbR = W_\bbR\ltimes \bbR^*_+$ the Langlands dual group,
$\bbS_\varphi$ be  the set of conjugacy classes of Langlands parameters modulo the connected component of identity map. For any $s\in 
\bbS_\varphi$,
$\check{H}_s = \Cent(s,\check G)^\circ$ the connected component of the centralizer of $s\in\bbS_\varphi$ we have $\check{H}_s$ is conjugate with $H$. 
Following the D. Shelstad pairing $$\langle s, \pi\rangle : \bbS_\varphi \times \Pi(\varphi) \to \bbC$$
$$\varepsilon(\pi) = c(s)\langle s,\pi\rangle.$$
Therefore, the relation
$$\sum_{\sigma\in \Sigma_s} \tr \sigma(f^H) = \sum_{\pi\in\Pi} \varepsilon(\pi) \tr \pi(f)$$
can be rewritten as
$$\widetilde{\Sigma}_s(f^H) = \sum_{s\in\Pi} \langle s,\pi\rangle \tr\pi(f)$$
and
$$\widetilde{\Sigma}_s(f^H) = c(s)^{-1}\sum_{\sigma\in{\Sigma}_s} \tr \sigma(f^H) .$$
We arrive, finally to the result
\begin{thm}\cite{labesse}
$$\tr \pi(f) = \frac{1}{\#\bbS_\varphi}\sum_{s\in\bbS_\varphi} \langle s,\pi\rangle \widetilde{\Sigma}_s(\check{f}).$$
\end{thm}

\begin{thm}
$$\tr R(f)|_{L^2_{cusp}(\Gamma \backslash \SU(2,1)} = \sum_{\Pi_\mu}\sum_{\pi\in\Pi_\mu}m(\pi)\mathcal S\Theta_\pi(f) = \sum_{ \Pi_\mu} \Delta(\gamma,\gamma_H)\mathcal {SO}(f_\mu),$$ where
$$\mathcal S\Theta_\pi(f) = \sum_{\pi\in \Pi}\kappa(\pi) \Theta_\pi(f)$$ is the sum of Harish-Chandra characters of the discrete series running over the stable conjugacy classes of $\pi$ 
and $$\mathcal {SO}(f_\mu) = \sum_{\lambda\in \Pi_\mu}\kappa(\pi_\lambda) \mathcal O(f_\lambda)$$ is the sum of orbital integrals weighted by a character $\kappa : \Pi_\mu \to \{\pm1\}$. 
\end{thm}
\beginpf\;
The proof just is  a combination of the previous theorems.
\endpf

\end{document}